\documentclass[amscd]{amsart}

\newtheorem{theorem}{Theorem}[section]
\newtheorem{lemma}[theorem]{Lemma}
\newtheorem{corollary}[theorem]{Corollary}
\newtheorem{definition}[theorem]{Definition}
\newtheorem{example}[theorem]{Example}
\newtheorem{remark}[theorem]{Remark}
\newtheorem{proposition}[theorem]{Proposition}

\usepackage{amsfonts}
\usepackage{amssymb}
\usepackage{amsmath}
\usepackage{eucal}
\usepackage{framed}
\usepackage[dvips]{graphicx}
\usepackage{epic}
\usepackage{graphics}
\usepackage{amsthm}%omake 
\usepackage{fancybox}%omake
\usepackage{color} 
\usepackage{ulem} 
\usepackage{eepic}

\begin{document}

\begin{flushright}
 {\tiny --[2026-OES=JC-essential+1FF-2-1.tex] --}
\end{flushright}

%\bigskip  
  
\title{On open embeddings of affine spaces in affine varieties and the Jacobian Conjecture $(JC_2)$}

\author{Susumu ODA}

\maketitle

 {\small
 \begin{center}
Professor Emeritus, Kochi University\\ 

{\small  334  Hiradani, Taki-cho, Taki-gun, Mie 519-2178, JAPAN\\
  ssmoda@ma.mctv.ne.jp or ssmoda427@gmail.com}
\end{center}
}

%%%%%%%%%%%%%%%%%%%%%%%%%   
%%%%%%%%%%%%%%%%%%%%%%%%%
% \bigskip
%{\small 
%\begin{center}
%(Dedicated to the Author's Father, Tadashi ODA,\\
% who passed away on August 25, 2010,\\
% and to Mother Kikue ODA,\\
% who passed away on July 6, 2026.)
%\end{center}
% }

 %\bigskip 
%%%%%%%%%%%%%%%%%%%%%%%%%%

 \begin{abstract} Let $k$ be a field of \uline{characteristic $0$}.  Our final goal is  the following faded problem~:\\ 
 {\sf The Jacobian Conjecture $(JC_2)$~:}
 {\it If  $f_1,  f_2$  are elements in a polynomial ring  $k[X_1,  X_2]$  over $k$ such that  $  \det(\partial f_i/ \partial X_j) $ is a nonzero constant, then  $k[f_1,  f_2] = k[X_1,  X_2]$.
}

  For this purpose, we consider the following Result~:
 
\noindent
{\sf Open Embeddings of Affine Spaces in Affine Varieties.}
{\it Let $X$ be a normal $k$-affine variety of dimension $2$ and let $U$ be  an open $k$-subvariety of $X$  which is isomorphic to $\mathbb{A}^2_k$.  Then $X = U$.}

This is effective for a more general result than our original objective $(JC_2)$.  

\noindent
{\sf The Generalized Jacobian Conjecture $(GJC_2)$.}   
{\it  Let $\varphi : X \rightarrow Y$ be an unramified morphism of normal $k$-affine varieties of dimension $2$. If both  $X$ and $Y$ are simply connected, then $\varphi$ is an isomorphism.}

 These results derived by $\mathbb{C}$-topological-method hold for  $\mathbb{C}$ instead of $k$ by ``Lefschetz-principle''. 
 \end{abstract}

%%%%%%%%%%%%%%%%%%%%%%%%%%%%%%%%%%%%

\renewcommand{\normalbaselines}{\baselineskip23pt \lineskip4pt \lineskiplimit5pt}
\newcommand{\longmapright}[1]{\smash{\mathop{\hbox to 1.1cm{\rightarrowfill}}\limits^{#1}}}
\newcommand{\mapright}[1]{ \xrightarrow{#1}}
\newcommand{\mapleft}[1]{ \xleftarrow{#1}}
\newcommand{\mapdown}[1]{ \Big\downarrow\llap{ $\vcenter{\hbox{$\scriptstyle#1\,$}}$ } }
\newcommand{\mapup}[1]{ \Big\uparrow\llap{ $\vcenter{\hbox{$\scriptstyle#1\,$}}$ } }
%\newcommand{\mapup}[1]{\xuparrow{#1}}

%%%%%%%%%%%%%%%%%%%%%%%%%%%%%%%%%%%%%
%%%%%%%%%%%%%%%%%%%%%%%%%%%%%%%%%%%%%%%%%%%%%%%%%%%%%%%%%%%%%%%%%%
\newcommand{\Liminj}[1]{\raisebox{-1.98ex}%
                        {$\stackrel{\normalsize\mbox{$\varinjlim$}}{\scriptstyle #1}$}}
\newcommand{\Limproj}[1]{\raisebox{-2.00ex}%
                        {$\stackrel{\normalsize\mbox{$\varprojlim$}}{\scriptstyle #1}$}}
                
%%%%%%%%%%%%%%%%%%%%%%%%%%%%%%%%%%%
%%%%%%%%%%%%%%%%%%%%%%%%%%%%%%%%%%%%%

\renewcommand{\thefootnote}{\fnsymbol{footnote}} %footnote

\footnote[0]{2020 {\it Mathematics Subject Classification}. Primary: 14R15; Secondary: 13A18 \\  
\ \ \ {\it Key words and phrases}. The Jacobian Conjecture, the Deep Jacobian Conjecture, Krull domains, \\
\ \ \ divisorial fractional ideals,  subintersections,  unramified, \'{e}tale, simply connected}

%%%%%%%%%%%%%%%%%%%%%%%%%%%%%%%%%%%% 
\markboth{\protect\footnotesize{{SUSUMU ODA}}}{\protect\footnotesize{OPEN EMBEDDINGS OF AFFINE SPACES AND THE JACOBIAN CONJECTURE }}   
%%%%%%%%%%%%%%%%%%%%%%%%%%%%%%%%%%%%

%\baselineskip=4.80mm 

\baselineskip=4.90mm

%%%%%%%%%%%%%%%%%%%%%%

\bigskip

{\bf Contents}
\begin{small}

\ \ \ 1. Introduction \dotfill \pageref{con2.1}

\ \ \ 2. Open Embeddings of $k$-Affine Spaces in  $k$-Affine Varieties 

 \ \ \ \ \ \ \  and The Jacobian Conjectures $(GJC_2)$ and $(JC_2)$ \dotfill  \pageref{con2.3}
 
\ \ \ 3. The Certain Extension of The Jacobian Conjecture$(JC_2)$ \dotfill \pageref{con2.4}

\ \ \ 4. Some Comments on The (Known) Example  \dotfill \pageref{conE.3} 

\ \ \ Appendix A.\ A Collection of Tools Required in This Paper \dotfill \pageref{con03}

\ \ \ Appendix B.\ A  (recent) Counterexample to the Jacobian Conjecture $(JC_3)$ \dotfill \pageref{con001}

\ \ \ References \dotfill \pageref{Bib}\\
\end{small}

%%%%%%%%%%%%%%%%%%%%%%%%%%%%%%%%%%%%%%

\vspace{2mm}  
              
{\bf {\large \section{Introduction}}} \label{con2.1}

 This paper is derived from investigating the Jacobian Conjecture and the Zariski Main Theorem. We begin to describe what is the Jacobian Conjecture. 
   
    Let $ k $ be an algebraically closed field,
let $ \mathbb{A}^n_k =  {\rm Spec^m}(k[X_1, \ldots, X_n])$ be an affine space of dimension $n$ over $ k $  and let  $ f : \mathbb{A}_k^n \longrightarrow \mathbb{A}_k^n $ be a morphism of affine spaces over $k$ of dimension $n$.  Note here that for a ring $R$, ${\rm Spec}(R)$ (resp. ${\rm Spec^m}(R)$) denotes the prime spectrum of $R$ (or merely the set of prime ideals of $R$)  (resp. the maximal spectrum (or merely  the set of the maximal ideals of $R$)). 
Then $ f $  is given by 
$$\mathbb{A}^n_k \ni(x_1, \ldots, x_n) \mapsto (f_1(x_1, \ldots,x_n), \ldots, f_n(x_1, \ldots,x_n)) \in \mathbb{A}^n_k,$$
where $ f_i(X_1, \ldots, X_n)  \in k[X_1,  \ldots, X_n].$
  If $ f $  has an inverse morphism,
then  the Jacobian $J(f) := \det( \partial f_i/ \partial X_j) $ is a nonzero constant.  This follows from the easy chain rule of differentiations without specifying the characteristic of $k$.  The Jacobian Conjecture asserts the converse.

If $ k $ is of characteristic \ $ p > 0 $  and $ f(X) = X  + X^p,$  then $ df/dX = f'(X) = 1 $ but
$ X $ can not be expressed as a polynomial in $ f$. It follows that the inclusion $k[X+X^p] \hookrightarrow k[X]$ is finite and \'{e}tale but $f : k[X] \rightarrow k[X]$ is not an isomorphism. 
 This implies that $k[X]$ is not simply connected ({\it i.e.,} ${\rm Spec}(k[X]) = \mathbb{A}^1_{\mathbb{C}}$ is not simply connected, see \S4.Definition \ref{def-3}) when ${\rm char}(k) = p > 0$. 
Thus we must assume that the characteristic of $ k $ is $0$.
 
\bigskip

The algebraic form of {\bf The Jacobian Conjecture$(JC_n)$} (or {\bf the Jacobian Problem}$(JC_n)$) is the following :
\bigskip

\noindent
{\bf The algebraic form $(JC_n)$.}  \begin{sl}
 If  $f_1,\ldots, f_n$  are elements in a polynomial ring  $k[X_1,\ldots,  X_n]$  over a field $k$ of characteristic $0$ such that  $\det(\partial f_i/ \partial X_j) $ is a nonzero constant, then  $k[f_1,\ldots, f_n] = k[X_1,\ldots,  X_n]$.
\end{sl}

\bigskip

Note that when considering $(JC_n)$, we may assume that $k=\mathbb{C}$ by ``Lefschetz-principle'' (See [10,(1.1.12)]).

%\bigskip

  The Jacobian Conjecture$(JC_n)$  has been settled affirmatively under a few special assumptions below (See [6]).   Let $k$ denote a field of characteristic $0$. We may assume that $k$ is algebraically closed.  Indeed, we can consider it in the case $k = \mathbb{C}$, the field of complex numbers. So we can use all  of the notion of Complex Analytic Geometry. But in this paper, we go forward with the algebraic arguments. 

 For example, under each of the following assumptions, the Jacobian Conjecture$(JC_n)$ has been settled affirmatively\ (Note that we may assume that $k$ is algebraically closed and of characteristic $0$)~:\\
{\bf Case(1)} $f : \mathbb{A}_k^n = {\rm Spec}(k[X_1, \ldots, X_n]) \rightarrow {\rm Spec}(k[f_1, \ldots, f_n]) = \mathbb{A}_k^n$ is injective~;\\
{\bf Case(2)} $ k(X_1,\ldots ,X_n) =  k(f_1,\ldots,f_n)$~;\\
{\bf Case(3)} $ k(X_1,\ldots ,X_n) $ is a Galois extension of $ k(f_1,\ldots,f_n) $~;\\
{\bf Case(4)}  $ \deg f_i \leq 2 $\  for all $ i $~;\\
{\bf Case(5)} $ k[X_1,\ldots ,X_n] $ is integral over $ k[f_1,\ldots,f_n]$.
 
\bigskip

 A fundamental reference for The Jacobian Conjecture $(JC_n)$ is [6] which includes the above Cases.

See also the reference [6] for a brief history of the developments and the state of the art again since it was first formulated and partially proved by Keller in 1939 ([13]), together with a discussion on several false proofs that have actually appeared in print, not to speak of so many other claims of prospective proofs being announced but proofs not seeing the light of the day. The  Jacobian Conjecture$(JC_n)$, due to the simplicity of its statement, has already fainted the reputation of leading to solution with ease, especially because an answer appears to be almost at hand, but nothing has been insight even for $n=2$.

 The conjecture obviously attracts the attention of one and all. It is no exaggeration to say that almost every makes an attempt at its solution, especially finding techniques from a lot of branches of mathematics such as algebra (Commutative Ring Theory), algebraic geometry/topology, analysis (real/complex) and so on, having been in whatever progress (big or small) that is made so far (cf. E. Formanek, Bass' Work on The Jacobian Conjecture, Contemporary Mathematics 243 (1999), 37-45).

For more recent arguments about The Jacobian Conjecture, we can refer to [W] and [K-M].

%%%%%%%%%%%%%%%%%%%%%%%%%%%%%%%%%%%%%%%%%%

\bigskip

 Throughout this paper,  unless otherwise specified, we use the following notations~:
\bigskip

\noindent
{\bf $\langle$ Basic Notations $\rangle$}

\noindent
\textbullet\  All fields, rings and algebras are assumed to be  commutative with unity.
 
--- For an integral domain $R$,

\noindent 
\textbullet\ A {\it factorial} domain is also called  a unique factorization domain, a UFD,\\
\textbullet\ An integral domain $R$ is called a {\it locally factorial} domain if for each maximal (prime) ideal $M$, $R_M$ is factorial,\\
\textbullet\  $R^{\times}$ denotes the set of units of $R$,\\
\textbullet\  $nil(R)$ denotes the {\it nilradical} of $R$, {\it i.e.,} the set of the nilpotent  elements of $R$,\\
\textbullet\  $K(R)$ denotes the total quotient ring (or the total ring of fractions) of $R$, that is, letting $S$ denote the set of all non-zerodivisors in $R$, $K(R):= S^{-1}R$,\\
\textbullet\  When $R$ is an integral domain, for $f\in R\setminus \{ 0 \}$\ $R_f := \{ r/f^n\ |\ r\in R, n \in \mathbb{Z}_{\geq 0} \}\ (\subseteq K(R))$, \\
%\textbullet\  Let  $A$ be a subring of $R$ and let $A_1$ and $A_2$ be $A$-subalgebras of $R$.   $A_1\cdot A_2$ denotes the image of $A_1\otimes_AA_2 \rightarrow R\ (a_1\otimes a_2 \mapsto a_1a_2)$, which is an $A$-subalgebra of $R$, \\
\textbullet\  ${\rm Ht}_1(R)$ denotes the set of all prime ideals of height one in $R$,\\
\textbullet\  ${\rm Spec}(R)$ denotes the {\it affine scheme} defined by $R$ (or merely the set of all prime ideals of $R$), and ${\rm Spec}^m(R)$ denotes the set of the maximal ideals of $R$,\\
\textbullet\  Let $A \rightarrow B$ be a ring-homomorphism and $p\in {\rm Spec}(A)$. Then  $B_p$ means $B\otimes_AA_p$. 

--- Let $k$ be an algebraically closed field.  

\noindent 
\textbullet\   A (separated) scheme over  $k$ is called a {\it $k$-(algebraic )scheme}. A $k$-scheme locally of finite type over $k$ is called   a  {\it algebraic variety} over $k$ or simply a {\it $k$-variety} if it is \uline{integral} ({\it i.e.,} irreducible and reduced).\\
\textbullet\  A $k$-variety $V$ is called a {\it $k$-affine variety} or an {\it affine variety over $k$} if it is $k$-isomorphic to an affine scheme ${\rm Spec}(R)$ for some \uline{$k$-affine domain} $R$ ({\it i.e.,} $R$ is a finitely generated domain over $k$). $R$ is called  the {\it coordinate ring} of $V$ and denoted by $K[V]$, and $K(V) := K(R)$.  In particular, a {\it $k$-affine space} $\mathbb{A}^n_k$ is ${\rm Spec}(k[X_1,\ldots,X_n])$ with a polynomial ring $k[X_1,\ldots,X_n]$. \\
\textbullet\  An integral, closed $k$-subvariety of codimension one in a $k$-variety $V$ is called  a {\it hypersurface} of $V$. \\
\textbullet\   A closed subset (possibly reducible) of pure codimension one in  a $k$-variety $V$ is called a {\it divisor} (or {\it a divisorial closed subset})  of $V$.   
In other words, a divisor of $V$ is a finite union of hypersurfaces of $V$.
 Thus an irreducible and reduced divisor ({\it i.e.,} a prime divisor) is the same as a hypersurface in our terminology.
 
%---  Let $A$ be a Krull domain. 
%
%\noindent  
%\textbullet\ $v_P(I):= \inf \{ v_P(a)\ |\ a\in I \}$, (which is non-zero for finitely many members $P \in {\rm Ht}_1(A)$ according to the finite character property of a defining family of a Krull domain),\\
%\textbullet\ ${\rm Supp}^*(I) := \{ P\in {\rm Ht}_1(A)\ |\ v_P(I) \not= 0 \}$, \  (For $Q \in {\rm Ht}_1(A)$, ${\rm Supp}^*(I) \not\ni Q \ \Leftrightarrow\ I_Q = A_Q$.)\\
%\textbullet\ for  fractional ideals $I_1, \ldots, I_n$ of $A$, $\widetilde{\prod}_{i=1}^nI_i$ denotes a product (ideal multiplication)  $I_1\cdots I_n\ (\subseteq K(A))$ of fractional ideals (not a direct product), which is also a fractional ideal of $A$. 

%%%%%%%%%%%%%%%%%%%%%%%%%%%%
\bigskip
 
In Section 2, we settle the following our main results by use of the usual $\mathbb{C}$-topological methods. 
 \uline{Note that when considering our main results, a field $k$ of characteristic $0$ can be replaced with $\mathbb{C}$ by ``Lefschetz-principle''} (See [10,(1.1.12)]).

 \bigskip
 
%\begin{corollary}[{\bf Open Embedding of Affine Space}]   \label{Aff}

\noindent
{\bf Open Embeddings of Affine Spaces in Affine Varieties.}
{\it Let $X$ be a normal $\mathbb{C}$-affine variety of dimension $2$ and let $U$ be  an open $\mathbb{C}$-subvariety of $X$  which  is isomorphic to $\mathbb{A}^2_{\mathbb{C}}$.  Then $X = U$.}\ {\bf (\ref{Aff})}
%\end{corollary}

 \bigskip
%%%%%

%\begin{theorem} [{\bf The Generalized Jacobian Conjecture $(GJC)$}]  \label{GJC} 

\noindent
{\bf The Generalized Jacobian Conjecture $(GJC_2)$.}
{\it Let $\varphi : X \rightarrow Y$ be an unramified morphism of normal $\mathbb{C}$-affine varieties of dimension $2$.  
  If both $X$ and $Y$ are simply connected, then $\varphi$ is an isomorphism.}\ {\bf (\ref{GJC})}
  
%\end{theorem}
\bigskip
 
%\begin{corollary}[{\bf The Jacobian Conjecture} $(JC_n)$] \label{JC}

\noindent
{\bf The Jacobian Conjecture $(JC_2)$.} 
 {\it If  $f_1, f_2$  are elements in a polynomial ring $\mathbb{C}[X_1, ,X_2]$ such that  $\det(\partial f_i/ \partial X_j) $ is a nonzero constant, then  $\mathbb{C}[f_1, f_2] = \mathbb{C}[X_1,  X_2]$.} \ {\bf (\ref{JC})}
%\end{corollary}

 \bigskip
 
 %%%%%%%%%%%%%%%%%%%%%%%%%%%%
  
By the way, the Jacobian Conjecture $(JC_n)$ is a problem concerning a polynomial ring over a field $k$ (characteristic $0$), so that investigating the structure of automorphisms ${\rm Aut}_k(k[X_1,\ldots,X_n])$ seems to be substantial. Any member of ${\rm Aut}_k(k[X_1,X_2])$ is known to be tame, but for $n\geq 3$ there exists a wild automorphism of $k[X_1,\ldots,X_n]$ (which was conjectured by M.Nagata with an explicit example and  was settled by Shestakov and Umirbaev(2003)).  

 In such a sense, to attain a positive solution of $(JC_2)$ by an abstract argument like this paper may be far from its significance  . . . . . . .

%\bigskip
%%%%%%%%%%%%%%%%%%%%%%%%%

%For the consistency of our discussion, we assert that the examples
%\footnote[2]{See {\sf arXive}:0706.1138v99[math. AC] 30 Nov 2022~: Some comments around the examples against the Deep Jacobian conjecture (with some revision)}
% appeared in the papers ([12], [2] and [20])  which  would be against our original target Conjecture$(DJC)$, are imperfect or incomplete counter-examples. 
%Concerning this, we  observe some comments about ``Example'' in {[10,(10.3) in p.305]} (See {\bf A BREAK 2}\ below).  
%% {[4,(10.3) in p.305]}

%%%%%%%%%%%%%%%%%%%%%%%%%

\bigskip

  Remark that we often say in this paper that

\begin{quote}
\noindent
{\sf  a ring $A$ is  ``simply connected'' if ${\rm Spec}(A)$ is simply connected, and a ring homomorphism $f : A \rightarrow B$ is ``unramified, \'{e}tale, an open immersion, a closed immersion, $\cdots\cdots$" when ``so" is its morphism  ${}^af : {\rm Spec}(B) \rightarrow {\rm Spec}(A)$, respectively.}
\end{quote}

%\bigskip

%%%%%%%%%%%%%%%%%%%%%%%%%%%%%%%%
%%%%%%%%%%%%%%%%%%%%%%%%%%%%%%%%

 \vspace{4mm} 

{\bf {\large \section{Open Embeddings of $k$-Affine Spaces in  $k$-Affine Varieties and The Jacobian Conjectures $(GJC_2)$ and $(JC_2)$}}} \label{con2.3}
 
%%%%%%  

Though we should be going without saying, we prepare some notations for our purpose.

\bigskip

\noindent
{\bf Notations~:} Let $k$ be an algebraically closed field. 
 For a  $k$-affine domain $R$ and $I$ its ideal, ${\rm Spec}^m(R)$ denotes the maximal-spectrum of $R$, and  $V^m(I)$ denotes $V(I) \cap {\rm Spec}^m(R) = \{ M \in {\rm Spec}^m(R)\ |\ I \subseteq M \}$.
 It is known that if $k[z_1,\ldots,z_n]$ is a polynomial ring, then the correspondence ${\rm Spec}^m(k[z_1,\ldots,z_n]) \ni M = (z_1-a_1, \ldots,z_n-a_n) \leftrightarrow (a_1,\ldots,a_n) \in k^n$ induces the isomorphism ${\rm Spec}^m(k[z_1,\ldots,z_n]) \cong k^n$ as  $k$-varieties, and for an ideal $J$ of $k[z_1,\ldots,z_n]$, $V^m(J) = V^m(\sqrt{J})$ corresponds to a (closed) algebraic set $\{ (a_1, \ldots,a_n) \in k^n\ |\ g(a_1, \ldots,a_n)=0\ (\forall g\in J) \}$ of $k^n = \mathbb{A}_k^n$\  (Hilbert's Nullstellensatz {\bf [10,(A.5.2)]}), which can be identified.

\bigskip

 Let $R$ be an integral domain with quotient field $K$ and let $I$ be an ideal of $R$.
 The set $S(I;R) := \{ f \in K\ |\ fI^n \subseteq R\ (\exists n \in \mathbb{Z}_{\geq 0}) \}$, which is an integral domain containing $R$.  
 For any integer $n \geq 0$, set $I^{-n} := \{ f \in K\ |\ fI^n \subseteq R \}$.   Then $S(I;R) = \bigcup_{n\geq 0}I^{-n}$.  We call $S(I;R)$ an {\it $I$-transform} of $R$, and abbreviate $S(I;R)$ to $S$ when there is no confusion.   We say that $S(I;R)$ is {\it finite} if $S(I;R) = R[I^{-n}]$ for some $n$. (See [16,Ch.V]).  

\bigskip
 
\begin{lemma}[{[16,Theorem 3',Ch.V]}] \label{NL}
 Let $k$ be an algebraically closed field.   Let $X$ be a $k$-affine variety defined by a $k$-affine domain $R$ and let $V$ be a closed set defined by an ideal $I$ of $R$. 
 Then the open subset $X\setminus V$ of $X$ is $k$-affine if and only if $1 \in IS$, where $S$ is the $I$-transform of $R$.
 In this case, $V$ is pure of codimension one and $S$ is finite, that is, the $k$-affine domain of $X\setminus V$.
\end{lemma}

\begin{lemma}[{[16,Theorem 5,Ch.V]}] \label{NL2}
 Let $k$ be an algebraically closed field.  Let $X$ be a normal $k$-affine variety of dimension $2$ and let $F$ be a divisor of  $X$. Then $X\setminus F$ is $k$-affine.
\end{lemma}

An ideal $J$ of a Krull domain\footnote[3]{Let $A$ be an integral domain which is contained in a field $K$. The integral domain $A$ is said to be {\it Krull domain} provided there is a family $\{ V_i \}_{i\in I}$ of principal valuation rings ({\it i.e.,} discrete rank one valuation rings), with $V_i \subseteq K$, such that 

{\rm (i)} $A = \bigcap_{i\in I}V_i$

{\rm (ii)} Given $0 \not= f \in A$, there is at most a finite number of $i$ in $I$ such that $f$ is not a unit in $V_i$.  

 Such a family as $\{ V_i \}_{i\in I}$ (or $\{ p_i\cap A\ |\ p_i\ \mbox{is\ a\ non-zero\ prime\ ideal\ of}\ V_i\ (i\in I) \}$) is called a {\it defining family} of a Krull domain $A$. The property {\rm (ii)} is called {\it the finite character } of $\{ V_i \}_{i\in I}$ (or $A$).  Note that $\{ A_P\ |\ P\in {\rm Ht}_1(A) \}$ (or ${\rm Ht}_1(A)$) is indeed a defining family of a  Krull domain $A$ (cf.[11,(1.9)]), where $A_P$ is called an {\it essential} valuation over-ring of $A$ $(\forall P\in {\rm Ht}_1(A))$. The defining family $\{ A_P\ |\ P\in {\rm Ht}_1(A) \}$ is the \uwave{minimal} one\ (cf.[14,(12.5)]), which is called the {\it essential defining family} of $A$.}
 $A$ is called {\it divisorial} if $A:_{K(A)}(A:_{K(A)}J) = J$ (See [11]). 
 We note that a Noetherian normal domain is a Krull domain.   
 
\begin{corollary}  \label{NL3}  
 Let $k$ be an algebraically closed field.  Let $X$ be a normal $k$-affine variety of dimension $2$ and let $F$ be a closed subset of $X$.  
 Then $X\setminus F$ is $k$-affine if and only if  $F$ is a divisor of $X$  if and only if a defining ideal $I$ of $F$ (that is, $F = V^m(I)$) is a divisorial ideal of $K[X]$.  In particular, for any $P\in {\rm Ht}_1(K[X])$ $X\setminus V^m(P)$ is an open $\mathbb{C}$-affine subvariety of $X$. 
\end{corollary}

The following fact is well-known.

\begin{lemma} [{[AAG,(1.2.10)]}] \label{NMM}
  Let $S$ be an affine scheme, and let  $(X,f)$ is a scheme separated over $S$. If $U_1$ and $U_2$ are affine open set of $X$, then $U_1\cap U_2$ is an affine open set of $X$. Moreover if $f:X\rightarrow S$, both $U_1 \rightarrow S$ and $U_2 \rightarrow S$ are morphisms of finite type,  then $U_1\cap U_2\rightarrow S$ is a morphism of finite type.
\end {lemma} 

\begin{proof} 
 Consider the diagonal $\Delta_{X/S}:X \rightarrow X\times_S X$, which is a closed immersion. Then $U_1\cap U_2 = \Delta_{X/S}^{-1}(\Delta_{X/S}(X)\cap (U_1\times_S U_2))$ and $\Delta_{X/S}(X)\cap (U_1\times_S U_2)$ is a closed subscheme of $U_1\times_S U_2$. 
Since $U_1\times_S U_2$ is an affine scheme,  $\Delta_{X/S}(X)\cap (U_1\times_S U_2)$ is an affine scheme.  Since $U_1\times_S U_2$ is an open subset of $X\times_sX$,  $\Delta_{X/S}(X)\cap (U_1\times_S U_2)$ is an open subset of $\Delta_{X/S}(X)$.  Thus $U_1 \cap U_2$ is isomorphic to $\Delta_{X/S}(X)\cap (U_1\times_S U_2)$ as schemes.  Hence $U_1\cap U_2$ is an affine scheme. 

 The latter follows from the consideration of morphisms of finite type (noting that a  closed immersion is of finite type). 
\end{proof}

\begin{corollary}  \label{NL4}  
 Let $k$ be an algebraically closed field.  Let $X$ be a normal $k$-affine variety of dimension $2$ and let $F_1$ and $F_2$ be   divisors of $X$. Then  $F_1\cup F_2$ is a divisor of $X$. 
\end{corollary}

\begin{proof} It is easy to see that  $(X\setminus F_1)\cap (X\setminus F_2) = X\setminus (F_1\cup F_2)$. 
 Since $X\setminus F_1$ and $X\setminus F_2$ are open $k$-affine varieties by Lemma \ref{NL2},  $(X\setminus F_1)\cap (X\setminus F_2)$ is an open $k$-affine variety by Lemma \ref{NMM} (Put $S={\rm Spec}(k)$).   So $F_1\cup F_2$ is a  divisor of $X$ by Corollary \ref{NL3}.
 \end{proof}

%%%%%%%%%%%%%%%%%%%%%%%

To make sure, we begin with the following definitions.

\begin{definition}[Unramified, \'{E}tale] \label{def-2}
{\rm 
 Let $f : A \rightarrow B$ be a ring-homomorphism \uline{of finite type} of Noetherian rings.
 Let $P \in {\rm Spec}(B)$ and put $P\cap A := f^{-1}(P)$, a prime ideal of $A$. 
 The homomorphism $f$ is called {\it unramified}\footnote[2]{In general, let $X$ and $Y$ be of locally Noetherian schemes and  let  $\psi : Y \rightarrow X$ be a morphism locally of finite type. If for $y \in Y$,  $\psi_y^* : \mathcal{O}_{X,\psi(y)} \rightarrow \mathcal{O}_{Y,y}$ is unramified at $y$, then $\psi$ is called {\it unramified} at $y \in Y$.   The set $R_{\psi} := \{ y\in Y\ |\ \psi_{y}^*\ {\rm  is\ ramified} \} \subseteq Y$ is called  {\it the ramification locus} of $\psi$ and $\psi(R_{\psi})\subseteq X$ is called the {\it the branch locus} of $\psi$. Note that the ramification locus $R_{\psi}$ defined here is often called  the branch locus of $\psi$ instead of $\psi(R_{\psi})$ in some texts\ (indeed, see e.g. [4] {\it etc.}).}
 at $P \in {\rm Spec}(B)$ if $PB_P = (P\cap A)B_P$  and $k(P) := B_P/PB_P$ is a finite separable field-extension of $k(P\cap A) := A_{P\cap A}/(P\cap A)A_{P \cap A}$. If $f$ is not unramified at $P$, we say $f$ is {\it ramified} at $P$. 
 The set $R_f:=\{ P\in {\rm Spec}(B)\ |\ {}^af\ {\rm is\ ramified\ at}\ P\in {\rm Spec}(B) \}$ is called the the {\it ramification locus} of $f$, which is a closed subset of ${\rm Spec}(B)$.   
 The homomorphism $f$ is called {\it \'{e}tale} at $P$ if $f$ is unramified and flat at $P$. The homomorphism $f$ is called {\it unramified} (resp. {\it \'{e}tale}) if $f$ is unramified (resp. \'{e}tale) at every $P\in {\rm Spec}(B)$.
 The morphism ${}^af :{\rm Spec}(B) \rightarrow {\rm Spec}(A)$ is called {\it unramified} (resp. {\it \'{e}tale}) if $f : A \rightarrow B$ is unramified (resp. \'{e}tale).
}
\end{definition}

\begin{definition}[(Scheme-theoretically or Algebraically) Simply Connected]  \label{def-3}
{\rm    
 A Noetherian ring $R$ is called {\it $($algebraically or scheme-theoretically$)$  simply connected}\footnote[3]{In general, let $X$ and $Y$ be  locally Noetherian schemes and let $\psi : Y \rightarrow X$ be a morphism locally of finite type. If $\psi$ is  finite and surjective, then  $\psi$ (or $Y$) is called a {\it (ramified) cover} of $X$\ (cf.[4,VI(3.8)]).  If a cover $\psi$ is \'{e}tale, $\psi$ is called an {\it \'{e}tale cover} of $X$. If every connected \'{e}tale cover of $X$ is isomorphic to $X$, $X$ is said to be {\it (scheme-theoretically or algebraically) simply connected}.  
Remark that if $X$ is an algebraic variety over $\mathbb{C}$ then `` $X$ is a (geometrically) simply connected in the usual $\mathbb{C}$-topology $\Rightarrow$ $X$ is (scheme-theoretically or algebraically) simply connected '', but in general  the converse `` $\Leftarrow$ '' does not hold. 
 If $X \subseteq \mathbb{A}^n_{\mathbb{C}}$, there is no need to distinguish between them.}
 if the following condition holds~: \ 
Provided  any `connected' ring $A$ ({\it i.e.,} ${\rm Spec}(A)$ is connected) with a finite  \'{e}tale ring-homomorphism $\varphi : R \rightarrow A$,   $\varphi$ is an isomorphism.
}
\end{definition}

To settle our main results, we consider them at the geometric view point.  
We use here the $\mathbb{C}$-topological argument. 
    
\bigskip

\noindent
{\bf NOTE.}~:  \uline{Note that when considering our main results, a field $k$ of characteristic $0$ can be replaced with $\mathbb{C}$ by ``Lefschetz-principle''} (See [10,(1.1.12)]).

\bigskip

 Note here that $\mathbb{P}^1_{\mathbb{C}}$ is decomposed into $\mathbb{A}^1_{\mathbb{C}} \cup \mathbb{A}^1_{\mathbb{C}}$ with the certain  glueing of two affine lines $\mathbb{A}^1_{\mathbb{C}}$ and that the projective line $\mathbb{P}^1_{\mathbb{C}}$ is simply connected but is not $\mathbb{C}$-affine indeed. 
 
\bigskip
 
%%%%%%%%%%%%%%%%%%%%%%%%%%%%%%%%%%%
 %%  In the part below, if $dim X \geq 3$ then the following istatements arr  false. %% 

Now we pose the following problems~:
\bigskip

\noindent
{\bf Problem I.}   Let $X$ be a normal $\mathbb{C}$-affine variety of dimension $n$ and let $X \hookrightarrow \mathbb{C}^m$ be  a closed immersion  \ $(\exists m\in \mathbb{Z}_{\geq  1})$. 
  Then for a divisor $F$ of $X$,  there exists a divisor  $F'$  of $\mathbb{C}^m$ with $F' \cap X = F$.
  
\bigskip

\noindent
{\bf Problem II.}  Let $X$ be a normal $\mathbb{C}$-affine variety of dimension $n$ and let $U$ be  an open $\mathbb{C}$-subvariety of $X$  which is isomorphic to $\mathbb{A}^n_{\mathbb{C}}$.  Then $X = U$.

\bigskip  

 If one of Problems were affirmative for each $n= \dim(X)$, the Jacobian Conjecture $(JC_n)$ could be proved, and two Problems could be equivalent to $(JC_n)$ (by consideration of  Zariski's Main Theorem (Lemma \ref{E})). 
   Unfortunately, it was shown that \uline{there exists a negative solution to Problems with $(JC_n)$ if $\dim (X) \geqq 3$} (by AI on July 2026) (which  will be discussed in Appendix B below), however if $\dim (X) \leq 2$, Problems and $(JC_2)$  have  affirmative solutions shown in the sequel: 

\begin{lemma}  \label{EMM}
  Let $X$ be a normal $\mathbb{C}$-affine variety of dimension $2$  
 and let $X \hookrightarrow \mathbb{C}^n$ be  a closed immersion  \ $(n\in \mathbb{Z}_{\geq 1})$. 
  Then for an irreducible  divisor $F$ of $X$,  there exists a divisor  $F'$  of $\mathbb{C}^n$ with $F' \cap X = F$.
\end{lemma}

\begin{proof}
 Let $K[X]$ denote  the coordinate ring of $X$  and let $F = V^m(I)$ for an ideal $I$ of $K[X]$.  
 Since $F$ is irreducible,  $I$ is a primary ideal of $K[X]$. 
  Let $\varphi: \mathbb{C}^{[n]}:=\mathbb{C}[Y_1,\ldots,Y_n] \rightarrow K[X]$ be a surjective homomorphism induced by the closed immersion $X \hookrightarrow \mathbb{C}^n$, let $G:= {\rm Ker}(\varphi)$ and let $I':=\varphi^{-1}(I)$, a primary ideal of $\mathbb{C}^{[n]}$.  Then
  $$ \mathbb{C}^{[n]}/G \cong K[X]\ \ \mbox{and}\ \ (\mathbb{C}^{[n]}/G)/(I'/G) \cong \mathbb{C}^{[n]}/I' \cong K[X]/I.$$
  Since $I'/G = (\overline{y_1},\ldots, \overline{y_r})(\mathbb{C}^{[n]}/G)$ for some $y_i \in I'\setminus G$, we have 
  $$I' = G + (y_1,\ldots, y_r)\mathbb{C}^{[n]},$$
 where we can assume that each $y_i$ is primary element in $\mathbb{C}^{[n]}$ because $I'$ is a primary ideal of $\mathbb{C}^{[n]}$.   
   Thus  
 $$V^m(I') = V^m(G) \cap V^m(\varphi^{-1}(I)) = V^m(G) \cap V^m(y_1\mathbb{C}^{[n]})\cap \ldots \cap V^m(y_r\mathbb{C}^{[n]})$$
 and 
  $$V^m(I) = V^m(\varphi(I')) = V^m(\varphi(y_1)K[X])\cap \ldots \cap V^m(\varphi(y_r)K[X]).$$
 %%%%%%%%%%%%%%
Since $X$ is normal,  every prime divisor of $ \overline{y_i}(\mathbb{C}^{[n]}/G)$ is of height one  by  Lemma \ref{S}  and  every prime divisor of $y_i\mathbb{C}^{[n]}$ is of height one   $(1 \leq i \leq r)$.
%%%%%%%%%%%%%%
  Thus $V^m(I') \subseteq V^m(y_i\mathbb{C}^{[n]})$ for all  $(1 \leq i \leq r)$, 
   so that $V^m(I) \subseteq V^m(\varphi(y_i)K[X])$. Since $F=V^m(I)$ is a divisor,   $\bigcup_{i=1}^rV^m(\varphi(y_i)K[X])$ is a  divisor by Corollary \ref{NL4}.  Since $F$ is irreducible by assumption, 
  we have 
$$F = \bigcap_{i=1}^rV^m(\varphi(y_i)K[X])= \bigcup_{i=1}^rV^m(\varphi(y_i)K[X])\  (= V^m(\varphi(y_j)K[X]\  (\exists j)).$$  
Then  putting $F' := \bigcup_{i=1}^rV^m(y_i\mathbb{C}^{[n]})$, $F'$ is a  divisor of $\mathbb{C}^n$  and $F'\cap X = F$.  Therefore the desired $F'$ exists.   

Therefore there exists a divisor  $F'$  of $\mathbb{C}^n$ with $F' \cap X = F$, and hence $(\mathbb{C}^n\setminus F')\cap X = X\setminus F$. 
  %%%%%%  
\end{proof}

\begin{corollary}  \label{SSC-1}
  Let $X$ be a normal $\mathbb{C}$-affine variety of dimension $2$ and let $F$ be an irreducible divisor of $X$. 
  Then  $X\setminus F$  is not simply connected.
\end{corollary}

\begin{proof}
  {\bf Suppose that $X\setminus F$ is simply connected.}
    We can use Lemma \ref{EMM}. 
  Since a loop  in $\pi_1(X\setminus F)$ is  a finite product of geometric generators\footnote[2]{Let $D \subseteq \mathbb{C}^n$ be a divisor, and let $y$ be a non-singular point of $D$.  Consider   a real plane $\Pi \subseteq \mathbb{C}^n$ intersecting $D$ transversely at $y$.  Let $C \subseteq \Pi$ be a circle of small radius with center at $y$.  It is well known that the fundamental group $\pi_1(\mathbb{C}^n\setminus D,o)$ is generated by loops $\gamma$ of the following form~: 
 $\gamma$ consists of a path $L$ joining the point $o$ with a point $y_1 \in C$, a loop around $y$ along $C$  beginning and ending at $y_1$, and  returning to $o$ along the path $L$ in the opposite direction. Such loops $\gamma$ (and the corresponding elements in $\pi_1(\mathbb{C}^n\setminus D)$) will be called \uline{geometric generators} ([12])}
 in $\pi_1(\mathbb{C}^n\setminus F')$, it is not trivial in $\pi_1(X\setminus F)$. Note that any element in $\pi_1(\mathbb{C}^n\setminus F')$ is finitely generated by some loops (geometric generators) around $F'$ (See [12,\S 3]).    
 It follows that $X$ has a loop around $F$, which is also a loop (geometric generator) around $F'$ in $\mathbb{C}^n$. This loop is not contractible in $\mathbb{C}^n\setminus F'$ by Corollary \ref{D-3}, so that it is not contractible in $X\setminus F$, {\bf a contradiction}. 
 Therefore $X\setminus F$ is not simply connected.
 \end{proof}

%%%%%%%%%%%%%%%%%%%%%%%%%%%%%%%%%%

\begin{proposition}   \label{SSC}  
 Let $X$ be a normal  $\mathbb{C}$-affine variety  of dimension $2$ and let $F$ be a divisor of $X$.   Then the open $\mathbb{C}$-subvariety $X\setminus F$ is $\mathbb{C}$-affine and  \uline{not} simply connected. 
\end{proposition}

\begin{proof}  Let $F = \bigcup_{i=1}^rF_i$ be an irreducible decomposition of  $F$ with irreducible components $F_i$.  Suppose that  for  $k$ ($1 \leq k \leq r$), $\bigcup_{i=1}^kF_i$ is a divisorial closed subset of $X$ and  $X\setminus (\bigcup_{i=1}^kF_i)$ is not simply connected.  
 Since an open immersion $X\setminus (\bigcup_{i=1}^{k+1}F_i) \hookrightarrow X\setminus (\bigcup_{i=1}^{k}F_i)$ induces a surjective homomorphism $\pi_1(X\setminus(\bigcup_{i=1}^{k+1}F_i)) \rightarrow \pi_1(X\setminus(\bigcup_{i=1}^{k}F_i))$ by the consideration of forgetting unnecessary loops around $F_{k+1}$ in $X$.     
Then  replacing  $X$ (resp. $F$) with  $X\setminus (\bigcup_{i=1}^kF_i)$  (resp. $\bigcup_{i=1}^{k+1}F_i$)  and $\pi_1(X\setminus F_1) \not=1$ (non-trivial) by Corollary  \ref{SSC-1} ,   it follows that  $X\setminus (\bigcup_{i=1}^{k+1}F_i) = (X\setminus (\bigcup_{i=1}^kF_i))\setminus F_{k+1}$ is $\mathbb{C}$-affine and  not simply connected by Corollaries \ref{SSC-1} and \ref{NL4}. Hence we come to the clonclusion that $X\setminus F$ is not simply connected by induction on $k$.
\end{proof}

%%%%%%%%%%%%%%%%%%%%%%%%%%%%%% 
 
\begin{corollary} \label{SC} 
  Let $X$ be a normal  $\mathbb{C}$-affine variety of dimension $2$, and let $U\ (\subsetneq X)$ be an open $\mathbb{C}$-affine subvariety of $X$. Then  $U$ is \uline{not} simply connected. 
 \end{corollary}
 
\begin{proof}
 Let $K[X]$ denote the coordinate ring of $X$, which is a $\mathbb{C}$-affine domain.  Then $F:= X\setminus U = V^m(I)$ for an ideal $I$ of $K[X]$.  Since $U = X\setminus F$ is $\mathbb{C}$-affine, every prime divisor of $I$ is of height $1$ by Lemma \ref{NL},  ({\it i.e.,} $I$ is  pure of codimension one). 
  So $F=V^m(I)$ is a divisor of an irreducible $\mathbb{C}$-affine variety $X$, and $U=X\setminus V^m(I)$ is not simply connected by Proposition \ref{SSC}.
  \end{proof}  

%%%%%%%%%%%%%%%%%%%%%%%%%%%%%%%%

%%%%  %%%%%

 \bigskip
 
\noindent
 {\bf NOTE 1~:} Let  $\mathbb{P}^n_{\mathbb{C}}$ denote the projective space and let $F$ be a hypersurface in  $\mathbb{P}^n_{\mathbb{C}}$.  Then $\mathbb{P}^n_{\mathbb{C}}\setminus F$ is simply connected if $F$ is a hyperplane, and is not simply connected if $F$  is a hypersurface  except a hyperplane (Corollary \ref{DD-3}). 
 
 \bigskip
 
 \noindent
 {\bf NOTE 2~:} Let $U \hookrightarrow V$ be an open immersion of normal
 $\mathbb{C}$-affine varieties. Then $V\setminus U$ is a a divisor. If $U$ is simply connected, then $V$ is also simply connected.  Thus  $K[V]^\times = K[U]^\times = \mathbb{C}^\times$ by Proposition \ref{1.3. Proposition}.
   
 (Note that the simple-connectivity of $U$ gives that of $V$.  Indeed, it is known that  $\pi_1(U,p_0) \rightarrow \pi_1(V,p_0)\ (p_0 \in U)$ is surjective by forgetting unnecessary loops around the  divisor $V\setminus U$.)
 
   \bigskip

As an easy corollary to Corollary \ref{SC}, we have the following~: 

\begin{theorem}[{\bf Open Embeddings of Affine Spaces in Affine Varieties}]   \label{Aff}
 Let $X$ be a normal $\mathbb{C}$-affine variety of ${\rm dim}(X) = 2$, and let $U$ be  an open $\mathbb{C}$-subvariety of $X$  which is isomorphic to $\mathbb{A}^2_{\mathbb{C}}$.  Then $X = U$.
\end{theorem}
                                                         
\begin{proof}  Let $F:= X\setminus U$.  
 {\bf Suppose that $F$ is not empty}.  Then $F$ is a non-empty divisorial closed subset (a divisor) of $X$ by Corollary \ref{NL3}.  Since $U = X\setminus F$,  $U$ is not simply connected by Corollary \ref{SC}. 
 However $U \cong \mathbb{C}^2$ is simply connected, {\bf a contradiction}.
\end{proof}

 %%%%%

\begin{theorem} [{\bf The Generalized Jacobian Conjecture $(GJC_2)$}]  \label{GJC} 
  Let $\varphi : X \rightarrow Y$ be an unramified morphism of normal $\mathbb{C}$-affine varieties  of  dimension $2$.
  If both $X$ and $Y$ are simply connected, then $\varphi$ is an isomorphism.
\end{theorem}

\begin{proof}  Let $T = K[X]$ and $S = K[Y]$. 
 Note first that $\varphi : S \rightarrow T$ is an \'{e}tale (and hence flat) homomorphism by Lemma \ref{A.2} and that $\varphi$ is injective by Lemmas \ref{FO} and \ref{I1}.  
 We can assume that $\varphi : S\rightarrow T$ is the inclusion $S \hookrightarrow T$.
 Let $C$ be the integral closure of $S$ in $K(T)$. Then $S\hookrightarrow C$ is finite  and $C$ is a  $\mathbb{C}$-affine domain by Lemma \ref{D}, and $C \hookrightarrow T$ is an open immersion by Lemma \ref{E}.  
 Thus $C=T$ by Corollary \ref{SC}, and $S\hookrightarrow T=C$ is \'{e}tale by Lemma \ref{A.2} and finite.  Hence $S=T$ because $S$ is (algebraically) simply connected.   
 \end{proof}
 
 On account of Remark \ref{2.4. Remark}, Theorem \ref{GJC} 
 resolves The Jacobian Conjecture $(JC_2)$~:

\begin{corollary}[{\bf The Jacobian Conjecture} $(JC_2)$] \label{JC}
 If  $f_1,  f_2$  are elements in a polynomial ring $k[X_1,X_2]$  over a field $k$ of characteristic $0$ such that  $\det(\partial f_i/ \partial X_j) $ is a nonzero constant, then  $k[f_1,  f_2] = k[X_1,  X_2]$.
\end{corollary}

%%%%%%%%%%%%%%%%%%%%%%%%%%%%%%%%%% 

Finally, related to Theorem \ref{Aff} above, we note that the following theorem ({\it i.e.,} Problem II can not hold) because $(JC_n)$ for $n \geq 3$ is false. (See  \S 3 or [UAI]). 

\begin{theorem}  \label{notOpen}  
 For each $n \in \mathbb{Z}_{n\geq 3}$, there exists an $n$-dimensional $\mathbb{C}$-affine variety $X\not\cong \mathbb{A}^n_{\mathbb{C}}$ such that $X$  contain an  open subvariety which is isomorphic to $\mathbb{A}^n_{\mathbb{C}}$. 
\end{theorem}

%%%%%%%%%%%%%%%%%%%%%%%%%%%%%%%%%%

\vspace{4mm}   

{\bf {\large \section{The Certain Extension of The Jacobian Conjecture $(JC_2)$}}}  \label{con2.4}
 
%\bigskip
 
In this section we enlarge a coefficient ring of a polynomial ring and consider the Jacobian Conjecture about it.  This is seen in [6,I(1.1)] by use of the observation on the formal inverse [6,III]. We can also see it in [10,(1.1.114)].   Our proof is simpler than that of [6] %%%%%
 even though considering only the case of integral domains.  
\bigskip

\begin{theorem}\label{2.3-. Theorem}
  Let  $A$  be an integral
 domain whose quotient field  $K(A)$  is of characteristic
 $0$.  Let $ f_1, \ldots ,f_n$  be elements of a
 polynomial ring $ A[X_1,\ldots,X_n]$  such that
$$f_i = X_i + ({\rm higher\ degree\ terms})\ \ \ \ \ \ (1 \leq i \leq n)\ \ \ \ \ (*).$$ 
If $K(A)[X_1,\ldots,X_n] = K(A)[f_1,\ldots,f_n]$, then $A[X_1,\ldots,X_n] = A[f_1,\ldots,f_n].$
\end{theorem}

\begin{proof}
  It  suffices to prove $ X_1, \ldots ,X_n \in A[f_1, \ldots ,f_n]$.

 We introduce
 a linear order in the set  $ \{ k:=(k_1,\ldots,k_n) | k_r \in \mathbb{Z}_{\geq 0}\ (1 \leq r \leq n)  \}$  of lattice points in ${\mathbb{R}_{\geq 0}}^n$ (where $\mathbb{R}$ denotes the field of real numbers)  in the following way :
 $$k=(k_1, \ldots ,k_n) > j=(j_1, \ldots ,j_n)\ {\rm if}\ k_r > j_r\ {\rm for\ the\ first\ index}\ r\ {\rm with}\ k_r \not= j_r.$$
 (This order is so-called the lexicographic order in ${\mathbb{Z}_{\geq 0}}^n$). 
\bigskip

\noindent
{\bf Claim.} Let $F(s):=\sum_{j=0}^s c_jf_1^{j_1}\cdots f_n^{j_n} \in A[X_1,\ldots,X_n]$ with $c_j \in K(A)$. Then $c_j \in A\ (0 \leq j \leq s.)$

(Proof.)  If $s=0 (=(0,\ldots,0))$, then $F(0)=c_0 \in A$.  

Suppose that for  $k (< s)$,  
 $c_j \in A\ (0 < j \leq k)$. Then $F(k) \in A[X_1,\ldots,X_n]$  by $(*)$, and 
 $F(s) - F(k) = G:=\sum_{j>k}^s c_j f_1^{j_1}\cdots f_n^{j_n} \in A[X_1,\ldots,X_n]$. 
 Let $k'=(k'_1,\ldots,k'_n)$ be the next member of $k$ ($k=(k_1,\ldots,k_n)<(k'_1,\ldots,k'_n)=k'$) with $c_{k'}\not= 0$.

 We must show $c_{k'} \in A$. 
 Note that $F(s) = F(k)+G$ with $F(k),\ G \in A[X_1,\ldots,X_n]$. 
 Developing $F(s):=\sum_{j=0}^s c_jf_1^{j_1}\cdots f_n^{j_n}\in A[X_1,\ldots,X_n]$ with respect to $X_1,\ldots,X_n$, though the monomial $X_1^{k'_1}\cdots X_n^{k'_n}$ with some coefficient in $A$ maybe appears in $F(k)$, it appears in only one place of $G$ with a coefficient $c_{k'}$ by the assumption $(*)$.  
 Hence the coefficient of the monomial  $X_1^{k'_1}\cdots X_n^{k'_n}$ in $F(s)$ is a form $b+c_{k'}$ with $b \in A$  because $F(k) \in A[X_1,\ldots,X_n]$.  Since $F(s) \in A[X_1,\ldots,X_n]$, we have  $b+c_{k'} \in A$ and hence $c_{k'} \in A$.  
 Therefore we have proved our Claim by induction.  
\bigskip

 Next, considering  $K(A)[X_1, \ldots, X_n] = K(A)[f_1, \ldots, f_n]$, we have  
$$  X_1 = \sum c_{j}f_1^{j_1} \cdots f_n^{j_n}$$ 
with  $c_{j} \in A$ by {\bf Claim} above.
 Consequently,  $X_1$  is in $A[f_1, \ldots ,f_n]$.   Similarly $X_2, \ldots ,X_n$  are in  $A[f_1,\ldots,f_n]$  and the assertion is proved completely.
 Therefore $A[f_1,\ldots,f_n] = A[X_1,\ldots,X_n]$. 
\end{proof}

%\bigskip

   The Jacobian Conjecture for $2$-variables  can be generalized as follows.

\begin{corollary}[{cf.[10,(1.1.18)]}] \label{2.3. Theorem}
  Let  $A$  be an integral domain whose quotient field  $K(A)$  is of characteristic
 $0$.  Let $ f_1, f_2$  be elements of a polynomial ring $ A[X_1,X_2]$  such that the Jacobian  $\det( \partial f_i/ \partial X_j)$ is in $A^\times$. Then $A[X_1,X_2] = A[f_1,f_2].$
\end{corollary}

\begin{proof}
  We see that $K(A)[X_1, X_2] = K(A)[f_1, f_2]$ by Corollary \ref{JC}. 
 It  suffices to prove $ X_1, X_2 \in A[f_1, f_2]$.
 We may assume that $f_i\ (i = 1, 2$ have no constant term. 
  Since $f_i \in A[f_1, f_2]$, 
  $$f_i = a_{i1}X_1 + a_{i2}X_2 + \mbox{(higher degree terms)}$$
 with $ a_{ij} \in A,$  where  $(a_{ij}) = (\partial f_i/\partial X_j)(0,\ldots ,0)$.
 The assumption implies that the determinant of the matrix  $(a_{ij})$  is a unit in $ A$.
 Let 
$$  Y_i = a_{i1}X_1 + a_{i2}X_2 \ \ \ (i = 1, 2).$$ 
Then  $A[X_1, X_2] = A[Y_1, Y_2]$
  and $ f_i = Y_i + \mbox{(higher degree terms)}$. 
 So to prove the assertion,  we can assume that without 
loss of generality 
$$f_i = X_i + ({\rm higher\ degree\ terms})\ (i =1, 2)\ \ \ \ \ \ \ \ \ \ (*).$$

%%%%%%%%%%%%%%
 Therefore by Theorem \ref{2.3-. Theorem} we have $A[f_1, f_2] = A[X_1, X_2]$. 
\end{proof}

\begin{example} \label{3.2.Corollary}  
  Let  $\varphi : \mathbb{A}^2_{\mathbb{Z}} \rightarrow  \mathbb{A}^2_{\mathbb{Z}}$ be a morphism of affine planes  over $\mathbb{Z}$, the ring of integers.
  If the Jacobian $J(\varphi)$ is equal to either \ $\pm 1$, then $\varphi$ is an isomorphism.
\end{example}

 %%%%%%%%%%%%%%%%%%%  

\vspace{4mm}

{\bf {\large \section{Some Comments on The (known) Example}}} \label{conE.3}

\bigskip

%%%%%%%%%%%%%%%%%%%%%%%%%%%%%%

---------- {\large \bf A BREAK }----------------- 
 
  We find  Example sited in [10,(10.3) in p.305] which is  concerned to an open embedding of the $\mathbb{C}$-affine space in a $\mathbb{C}$-affine variety ([11,p.32]).  
 We investigate it in this section. 

\bigskip
 
{\bf Example} (An open embedding of $\mathbb{C}^2$ in a $2$-dimensional
 $\mathbb{C}$-affine variety).   

  Let $X := \{ (s,t,u)\in \mathbb{C}^3\ |\ su-t^2+t = 0 \}$, a closed  algebraic set in $\mathbb{A}^3_{\mathbb{C}} = \mathbb{C}^3$ 

 and 
 let $F :  \mathbb{C}^2 \rightarrow X$ be  a polynomial map given by 
 $$ F(x,y) := (y, xy, x^2y-x).$$
\ \ \ \  Then $F$ is an open embedding ({\it open immersion}) of $\mathbb{C}^2$ in ({\it into}) $X$.   

\ More precisely
 $$(F~:)\   \mathbb{C}^2 \cong \{ (s,t,u) \in X\ |\ s \not= 0 \} \cup \{ (s,t,u) \in  X\ |\ t-1\not= 0 \}.$$

\hfill -------------------\ \ \ \ \ \  

\bigskip   
  
\noindent
{\bf Note.}  A word `embedding' is defined in [10,(5.3.1)] and a word `open immersion' is defined in [10,p.285],  but  we can not find a word `{\it open embedding}'. 

In {\bf Example}, the parts $(\ \ \ )$ are the author's interpretations. 

This example fascinates us, however it could be possibly a counter-example to  Theorem \ref{Aff}).  
 
\bigskip 

To make sure, \uline{we will check Example above in more detail}.

{\bf Suppose that Example above is valid}.

It is easy to see that $X$ is a non-singular $\mathbb{C}$-affine variety by the Jacobian criterion.  
 So its ring $\mathbb{C}[\bar{s},\bar{t},\bar{u}]$ of regular functions on $X$ is a regular domain and hence is a locally factorial domain, 
 where $\bar{s}, \bar{t}, \bar{u}$ is  the images of $s, t, u$ by the canonical homomorphism $\mathbb{C}[s,t,u] \rightarrow \mathbb{C}[s,t,u]/(su-t^2+t) (= \mathbb{C}[\bar{s},\bar{t},\bar{u}])$.  %%%%%
 (Thus any $P\in {\rm Ht}_1(\mathbb{C}[\bar{s},\bar{t},\bar{u}])$ is an invertible ideal of $\mathbb{C}[\bar{s},\bar{t},\bar{u}]$ by {[11,(9.2)]}.)

The morphism $F$  induces $F^* : \mathbb{C}[\bar{s},\bar{t},\bar{u}] \rightarrow \mathbb{C}[x,y]$ by $F^*(\bar{s}) = y, F^*(\bar{t}) = xy$ and $F^*(\bar{u}) = x^2y-x$. 
 Then  $F^* : \mathbb{C}[\bar{s},\bar{t},\bar{u}] \overset{F^*}{\longrightarrow} \mathbb{C}[y,xy,x^2y-x] \hookrightarrow \mathbb{C}[x,y]$. %%%%%%%%
 
We see that $\{ (s,t,u) \in X\ |\ s\not=0 \}$ and $\{ (s,t,u) \in X\ |\ t\not=1 \}$ are open $\mathbb{C}$-affine subvarieties of $X$ and hence that $\{ (s,t,u) \in X\ |\ s\not=0 \}\cup \{ (s,t,u) \in X\ |\ t\not=1 \}$ is an open subset of $X$. 
 Since we supposed that {\bf Example} is valid, $F~: {\rm Spec}^m(\mathbb{C}[x,y]) = \mathbb{C}^2 \cong \{ (s,t,u) \in X\ |\ s \not= 0 \} \cup \{ (s,t,u) \in  X\ |\ t-1\not= 0 \}$. 
  
 Note that  $\{ (s,t,u) \in X\ |\ s\not=0 \} = {\rm Spec}^m(\mathbb{C}[s,t,u]_s)\cap X = {\rm Spec}^m(\mathbb{C}[\bar{s},\bar{t},\bar{u}]_{\bar{s}})$ and $\{ (s,t,u) \in X\ |\ t\not=1 \} = {\rm Spec}^m(\mathbb{C}[s,t,u]_{t-1})\cap X = {\rm Spec}^m(\mathbb{C}[\bar{s},\bar{t},\bar{u}]_{\bar{t}-1})$.

So we see the following~:\\
\ \ \ \  $ \{ (s,t,u) \in X\ |\ s\not=0 \} \cup \{ (s,t,u) \in X\ |\ t\not=1 \}$  \\
\ \ \ \ $= {\rm Spec}^m(\mathbb{C}[\bar{s},\bar{t},\bar{u}]_{\bar{s}}) \cup {\rm Spec}^m(\mathbb{C}[\bar{s},\bar{t},\bar{u}]_{\bar{t}-1})$\\
\ \ \ \ $= {\rm Spec}^m(\mathbb{C}[\bar{s},\bar{t},\bar{u}])\setminus (V^m(\bar{s})\cap V^m(\bar{t}-1))$\\
\ \ \ \ $= {\rm Spec}^m(\mathbb{C}[\bar{s},\bar{t},\bar{u}]) \setminus V^m(\bar{s},\bar{t}-1)$,\\
 which is an open subvariety of ${\rm Spec}^m(\mathbb{C}[\bar{s},\bar{t},\bar{u}]) = X$.
%%%

 The ideal $(\bar{s},\bar{t}-1)\mathbb{C}[\bar{s},\bar{t},\bar{u}]$ is a prime ideal of height $1$ and $V^m(\bar{s},\bar{t}-1)$ is isomorphic to a line $\mathbb{A}^1_{\mathbb{C}}$, a contractible hypersurface of $X$ in the usual $\mathbb{C}$-topology.

\noindent   
 [Indeed,   $\mathbb{C}[\bar{s},\bar{t},\bar{u}]/(\bar{s},\bar{t}-1) = \Bigl(\mathbb{C}[s,t,u]/(su-t^2+t)\Bigr)/\Bigl((s,t-1)/(su-t^2+t)\Bigr) \cong \mathbb{C}[s,t,u]/(s,t-1) \cong \mathbb{C}[u]$ and $V^m(\bar{s},\bar{t}-1)$ is of codimension $1$ in $X$.] 

 Thus the prime ideal $(\bar{s},\bar{t}-1)$ is in ${\rm Ht}_1(\mathbb{C}[\bar{s},\bar{t},\bar{u}])$ and hence is  a divisorial ideal of the  regular  domain  $\mathbb{C}[\bar{s},\bar{t},\bar{u}]$. So it is an invertible ideal.   
 Similarly, the ideal $(\bar{s},\bar{t})\mathbb{C}[\bar{s},\bar{t},\bar{u}]$ is a prime ideal of height $1$ and $V^m(\bar{s},\bar{t})$ is isomorphic to a line $\mathbb{A}^1_{\mathbb{C}}$, a contractible hypersurface of $X$ in the usual $\mathbb{C}$-topology.

\bigskip

It is easy to see that $F$ is a non-surjective open immersion because 
$$F({\rm Spec}^m(\mathbb{C}[x,y])) \not\ni (0,1,c) \in X\ \ (\forall c\in \mathbb{C}).$$  

We may identify $\bar{s}, \bar{t}, \bar{u}$ with $F^*(\bar{s}), F^*(\bar{t}), F^*(\bar{u}) \in \mathbb{C}[x,y]$, respectively   and  $F^* : \mathbb{C}[\bar{s},\bar{t},\bar{u}] = \mathbb{C}[y,xy,x(xy-1)] \hookrightarrow \mathbb{C}[x,y]$. 
 So $\bar{s} = y, \bar{t} = xy, \bar{u} = x^2y-x$.

\bigskip
 
Incidentally, a $\mathbb{C}$-automorphism $\sigma$ of  $\mathbb{C}[s,t,u]$ defined by $\sigma(s)=s,\ \sigma(t) = 1-t,\ \sigma(u) = u$  induces a $\mathbb{C}$-automorphism of $\mathbb{C}[\bar{s},\bar{t},\bar{u}]$, where we use the same $\sigma$. 
[Indeed, $\sigma(su-t(t-1)) = su - (1-t)\bigl((1-t)-1)\bigr) = su - (1-t)(-t) = su-t(t-1)$.]

 Then ${}^a\sigma \in {\rm Aut}({\rm Spec}^m(\mathbb{C}[\bar{s},\bar{t},\bar{u}]))$ with ${}^a\sigma^2 = id_X$. 
Besides, $\sigma$ can be seen an automorphism of the quotient field $\mathbb{C}(\bar{s},\bar{t},\bar{u}) = \mathbb{C}(x,y)$.
 
  We see
 $${\rm Spec}^m(\mathbb{C}[x,y])) \overset{F}{\cong} X\setminus V^m(\bar{s},\bar{t}-1) \overset{{}^a\sigma}{\cong} X\setminus {}^a\sigma(V^m(\bar{s},\bar{t}-1)) = X\setminus V^m(\bar{s},\bar{t})$$
and  
$${}^a\sigma(V^m(\bar{s},\bar{t}-1)) = V^m(\bar{s},\bar{t}).$$
  Since $(\bar{s},\bar{t})\mathbb{C}[\bar{s},\bar{t},\bar{u}] + (\bar{s},\bar{t}-1)\mathbb{C}[\bar{s},\bar{t},\bar{u}] = \mathbb{C}[\bar{s},\bar{t},\bar{u}]$, it follows that $V^m(\bar{s},\bar{t}) \cap V^m(\bar{s},\bar{t}-1) = \emptyset$ and $\bigl(X\setminus V^m(\bar{s},\bar{t})\bigr) \cup \bigl(X\setminus V^m(\bar{s},\bar{t}-1)\bigr) = X\setminus (V^m(\bar{s},\bar{t}-1)\cap V^m(\bar{s},\bar{t})) = X$. 
 Thus  
 $$ X = F({\rm Spec}^m(\mathbb{C}[x,y]))\cup {}^a\sigma F({\rm Spec}^m(\mathbb{C}[x,y])),$$
 where ${\rm Spec}^m(\mathbb{C}[x,y]) \cong \mathbb{C}^2 \cong {}^a\sigma F({\rm Spec}^m(\mathbb{C}[x,y]))$.

 Therefore  $X =  F({\rm Spec}^m(\mathbb{C}[x,y]))$ by Theorem \ref{Aff}.  However $X \not=  F({\rm Spec}^m(\mathbb{C}[x,y]))$ as was seen before. This is {\bf a contradiction.}

%%%%%%%%%%%%%%%%%%%%%%%%%%%%%%%%%

\vspace{2mm}

\appendix
 
\vspace{4mm}

{\bf {\large \section{A Collection of Tools Required in This Paper}}}  \label{con03}

\bigskip   

%%%%%%%%%%%%%%%%%%%%
 Recall the following well-known results, which are required in this paper.  We write down them for convenience.

\begin{remark}[{cf.[10,(1.1.31)]}] \label{2.4. Remark} 
{\rm   
    Let  $k$  be an algebraically closed  field of characteristic $0$ and let  $k[X_1, \ldots, X_n]$ denote a polynomial ring  and let  $f_1, \ldots, f_n \in k[X_1, \ldots, X_n]$.  If  the Jacobian $\det(\partial f_i/\partial X_i) \in k^\times ( = k \setminus (0))$, then $k[X_1, \ldots, X_n]$ is \'{e}tale over the subring $k[f_1, \ldots, f_n]$.
 Consequently $f_1, \ldots, f_n$  are algebraically independent over $k$.  
  
 In fact, put  $T = k[X_1, \ldots, X_n]$  and $S = k[f_1, \ldots, f_n] ( \subseteq T)$.
We have an exact sequence by [15,(26.H)]:
$$ \Omega^1_{k}(S)\otimes_ST \mapright{v} \Omega^1_{k}(T) \mapright{} \Omega^1_S(T) \mapright{} 0,$$
where 
 $$v(df_i\otimes 1) = \sum_{j=1}^n\dfrac{\partial f_i}{\partial X_j}dX_j\ \ \ \ \ (1 \leq i \leq n).$$
  So $\det(\partial f_i/\partial X_j) \in k^\times$ implies that  $v$ is an isomorphism.  Thus  $\Omega^1_S(T) = 0$  and hence  $T$  is unramified over  $S$ by [4,VI,(3.3)]. So $T$ is \'{e}tale over $S$ by Lemma \ref{A.2} below. 
 Thus $df_1,\ldots,df_n \in \Omega^1_k(S)$ compounds a free basis of $T\otimes_S\Omega^1_k(S) = \Omega^1_k(T)$, which means that $K(T)$  is algebraic over  $K(S)$ and  that $f_1, \ldots, f_n$  are algebraically independent over  $k$. 
 }
\end{remark}

%%%%%%%%%%%%%%%%%%%%%%%%

 The following proposition is related to the  `{\it simple-connectivity}' of   affine spaces $\mathbb{A}^n_k\ (n \in \mathbb{Z}_{\geq 0})$ over a field $k$ of characteristic $0$. Its (algebraic) proof is given  without the use of the geometric fundamental group $\pi_1(\ \ )$ after embedding $k$ in $\mathbb{C}$ (the Lefschetz Principle). 

\begin{proposition}[{[23]}] \label{1.2. Proposition}  Let  $k$  be an algebraically closed field of characteristic $0$. Then  a polynomial ring  $k[Y_1, \ldots, Y_n]$ over $k$ is (algebraically) simply connected.
 \end{proposition}

\begin{proposition}[{[2,Theorem 3]}] \label{1.3. Proposition} 
 Any invertible regular function on a normal, (algebraically) simply connected  $\mathbb{C}$-variety is constant.
\end{proposition}

The following is well-known, but we write it down here for convenience.

\begin{lemma}[{[14]}]   \label{D} 
  Let $k$ be a field, let  $R$  be a $k$-affine domain and let  $L$ be a finite  algebraic field-extension of $K(R)$. 
Then the integral closure $R_L$ of $R$ in $L$ is finite over  $R$. 
\end{lemma}

 Moreover the above lemma can be generalized as follows.

\begin{lemma}[{[15,(31.B)]}] \label{BB}
 Let $A$ be a Noetherian normal domain with quotient field $K$, let $L$ be a finite separable algebraic extension field of $K$ and let $A_L$ denote the integral closure of $A$ in $L$.  Then $A_L$ is  finite over $A$.
\end{lemma}

%%%%%%%%%%%%%%%%%%%%

\begin{lemma}[The Approximation Theorem for Krull Domains\ {[11,(5.8)]}] \label{AppKrull}
 Let $A$ be a Krull domain.  Let $n(P)$ be a given integer for each $P$ in ${\rm Ht}_1(A)$ such that $n(P)= 0$ for almost all $P$. For any preassigned set $P_1,\ldots, P_r$ there exists $x\in K(R)^\times$ such that $v_P(x) = n(P_i)$ with $v_P(x)\geq 0$ otherwise, where $v_P(\ \ )$ denotes the (additive) valuation associated to the DVR $A_P$.
\end{lemma}

\begin{lemma}[{[15,(6.D)]}] \label{I1} 
 Let $\varphi : A \rightarrow B$ be a homomorphism of rings. Then ${}^a\varphi : {\rm Spec}(B) \rightarrow {\rm Spec}(A)$ is dominating (or dominant) $(${\it i.e.,} ${}^a\varphi({\rm Spec}(B))$ is dense in ${\rm Spec}(A))$ if and only if $\varphi$ has a kernel $\subseteq {\rm nil}(A):= \sqrt{(0)_A}$.
 If, in particular, $A$ is reduced, then  ${}^a\varphi$ is dominating $\Leftrightarrow$ ${}^a\varphi({\rm Spec}(B))$ is dense in ${\rm Spec}(A)$ $\Leftrightarrow$ $\varphi$ is injective. 
\end{lemma}

\begin{lemma}[{[14,(9.5)], [15,(6.I)]}]  \label{FO}
  Let $A$ be a Noetherian ring and let $B$ be an $A$-algebra of finite type. If $B$ is flat over $A$, then the canonical morphism $f : {\rm Spec}(B) \rightarrow {\rm Spec}(A)$ is an open map.
 (In particular, if $A$ is reduced ({\it eg.,} normal) in addition, then $A\rightarrow B$ is injective.)
\end{lemma}

%%%%%%%%%%%%%%%%%%%%

For a Noetherian ring $R$, the definitions of its {\it normality} (resp. its {\it regularity}) is  seen in [15,p.116], that is, $R$ is {\it a normal ring} (resp. {\it a regular ring}) if $R_p$ is a normal domain (resp. a regular local ring) for every $p\in {\rm Spec}(R)$.
    
\begin{lemma}[{[14,(23.8)], [15,(17.I)]} (Serre's Criterion on normality)]  \label{S} 
 Let $A$ be a Noetherian ring.  Consider the following conditions~:\\
$(R_1)$~: $A_p$ is regular for all $p \in {\rm Spec}(A)$  with ${\rm ht}(p) \leq 1$ ;\\
$(S_2)$~: ${\rm depth}(A_p)\geq {\rm min}({\rm ht}(p), 2)$  for all $p \in {\rm Spec}(A)$.

Then  $A$ is a normal ring if and only if $A$ satisfies $(R_1)$ and $(S_2)$.
\noindent
$($ Note that $(S_2)$ is equivalent to the condition that any prime divisor of $fA$ for any non-zerodivisor $f$ of $A$ is not an embedded prime.$)$ 
\end{lemma}
     
%%%%%%%%%%%%%%%%%%%%%%

%\begin{lemma}[cf.{[14,(23.9)]}]   \label{SS}
% Let $(A,m)$ and $(B,n)$ be Noetherian local rings and $A\rightarrow B$ a local homomorphism.  Suppose that $B$ is flat over $A$. 
% Then \\
%{\rm (i)} if $B$ is normal (or reduced), then so is $A$,\\
%{\rm (ii)} if both $A$ and the fiber rings of $A\rightarrow B$ are normal (or reduced), then so is $B$.
%\end{lemma}

 %%%%%%%%%%%%%%%%%%%%%%

%%%%%%%%%%%%%%%%%%%%%%
 
\begin{lemma}[{[SGA,(Expos\'{e} I, Cor.9.11)]}] \label{A.2}  
 Let $S$ be a Noetherian normal domain, let $R$ is an integral domain   and let $\varphi : S \rightarrow R$ be a ring-homomorphism of finite type.  If $\varphi$ is unramified, then $\varphi$ is \'{e}tale.
\end{lemma}

\begin{lemma}[{[21,p.42]} (Zariski's Main Theorem)] \label{E}  
  Let  $A$ be a ring and let $B$ be an $A$-algebra of finite type which is quasi-finite over  $A$.   Let  $\overline{A}$  be the integral closure of  $A$  in $B$.  Then  the canonical morphism  ${\rm Spec}(B) \rightarrow {\rm Spec}(\overline{A})$  is an open immersion.
 \end{lemma}

\begin{lemma}[{[9,Prop(4.1.1)]}] \label{D-1}
 Let $W$ be a  (possibly, reducible) quasi-projective subvariety of $\mathbb{P}^n_{\mathbb{C}}$ and let $\overline{W}$ be its closure. Then the following hold~:\\
{\rm (i)} $\pi_1(\mathbb{P}^n_{\mathbb{C}}\setminus W) = 0$ if $\dim(W)<n-1$~;\\
{\rm (ii)} $\pi_1(\mathbb{P}^n_{\mathbb{C}}\setminus W) = \pi_1(\mathbb{P}^n_{\mathbb{C}}\setminus \overline{W})$ if $\dim(W)=n-1$.
\end{lemma}

\begin{lemma}[{[9,Prop(4.1.3)]}] \label{D-2}
 Let $V_i\ (1 \leq i \leq k)$ be different  divisor of $\mathbb{P}^n_{\mathbb{C}}$ which have   $\deg(V_i)=d_i$. Let $V:= \bigcup_{i=1}^kV_i$.   Then 
 $$\pi_1(\mathbb{P}^n_{\mathbb{C}}\setminus V)/[\pi_1(\mathbb{P}^n_{\mathbb{C}}\setminus V), \pi_1(\mathbb{P}^n_{\mathbb{C}}\setminus V)] = H_1(\mathbb{P}^n_{\mathbb{C}}\setminus V) = \mathbb{Z}^{k-1}\oplus (\mathbb{Z}/(d_1,\ldots, d_k)\mathbb{Z}),$$
where $(d_1,\ldots, d_k)$ denotes the greatest common divisor and $[\  , \ ]$ denotes a  commutator group.  
\end{lemma}

\begin{corollary}[{[9,Prop(4.1.4)]}] \label{D-3}
 If $V \subseteq \mathbb{C}^n$ is a divisor  with $r$ irreducible components, then 
$$ \pi_1(\mathbb{C}^n\setminus V) \rightarrow \mathbb{Z}^r$$
 is surjective.
\end{corollary}

\begin{corollary} \label{DD-3}
  Let $V_i\ (1 \leq i \leq k)$ be different  divisor of $\mathbb{P}^n_{\mathbb{C}}$ which have   $\deg(V_i)=d_i$. Let $V:= \bigcup_{i=1}^kV_i$. 
 Then   $\mathbb{P}^n_{\mathbb{C}}\setminus V$ is simply connected $\Longleftrightarrow$ $V$ is a hyperplane in $\mathbb{P}^n_{\mathbb{C}}$ $\Longleftrightarrow$ $\mathbb{P}^n_{\mathbb{C}}\setminus V \cong \mathbb{A}^n_{\mathbb{C}}$.
\end{corollary}

\begin{proof} 
 By Lemma \ref{D-2}, $\mathbb{P}^n_{\mathbb{C}}\setminus V$ is simply connected if and only if  $k=1$ and $d_1=\deg(V)=1$ if and only if  $V$ is a hyperplane in $\mathbb{P}^n_{\mathbb{C}}$  if and only if $\mathbb{P}^n_{\mathbb{C}}\setminus V \cong \mathbb{A}^n_{\mathbb{C}}$.
 \end{proof}

\vspace{4mm}

{\bf {\large \section{ A (recent) Counterexample to the Jacobian Conjecture $(JC_3)$}}}  \label{con001}

\bigskip   

%%%%%%%%%%%%%%%%%%%%

On July 20, 2026, L. Alp\"{o}ge at Anthropic (with Claude Fable 5) ([AI])  says on $\mathbb{X}$  that Example gives a polynomial map $(f_1,f_2,f_3) : \mathbb{C}^3 \rightarrow \mathbb{C}^3$ which is \'{e}tale (because its Jacobian determinant is a constant $-2$) but not injective ($\because [\mathbb{C}(x_1,x_2,x_3) : \mathbb{C}(f_1,f_2,f_3)] = 3$). So it is a Keller map which is not injective.

\bigskip

\begin{example} [{[AI]}, {[UAI]}]  \label{Ex}    
 {\rm 
 Let $F=(f_1,f_2,f_3) : \mathbb{C}^3 \rightarrow \mathbb{C}^3$ be a polynomial map defined by 
\begin{eqnarray*}
f_1[x_1,x_2,x_3] &:=& (1+x_1 x_2)^3 x_3+x_2^2 (1+x_1 x_2) (4+3x_1x_2)\\
&=&4 x_2^2 + 7 x_1 x_2^3 + 3 x_1^2 x_2^4 + x_3 + 3 x_1 x_2 x_3 + 
    3 x_1^2 x_2^2 x_3 + x_1^3 x_2^3 x_3,\\
         &   &       \\
f_2[x_1,x_2,x_3] &:=& x_2+3x_1(1+x_1x_2)^2x_3+3x_1x_2^2(4+3x_1x_2)\\
&=&x_2 + 12 x_1 x_2^2 + 9 x_1^2 x_2^3 + 3 x_1 x_3 + 6 x_1^2 x_2 x_3 + 
    3 x_1^3 x_2^2 x_3,\\
   &   &   \\
f_3[x_1,x_2,x_3] &:=& 2x_1-3x_1^2x_2-x_1^3x_3.
\end{eqnarray*}
 Then ${\rm det Jac}(F) = -2$ and $F(0,0,-1/4) = F(1,-3/2,13/2) = F(-1,3/2,13/2) = (-1/4,0,0)$.
}
\end{example} 

\bigskip
Put $A = 1+x_1x_2$ and $s=x_1/A$.  Then $f_2 = x_2+3f_1s, \ \ \ f_3 = 2s-x_2s^2 - f_1s^3$.
 On the dense open set $A\not= 0$,  the coordinate changes : 
$$ (x_1,x_2,x_3) \mapsto (f_1,x_2,s) \mapsto (f_1,f_2,f_3)$$  
have determinants
$$ {\rm det}\dfrac{\partial(f_1,x_2,s)}{\partial(x_1,x_2,x_3)} =-A, \ \ \ \ {\rm det}\dfrac{\partial(f_1,f_2,f_3)}{\partial(f_1,x_2,s)}
=2/A.$$
 Thus  ${\rm det}\dfrac{\partial(f_1,f_2,f_3)}{\partial(x_1,x_2,x_3)} = {\rm det}\dfrac{\partial(f_1,f_2,f_3)}{\partial(f_1,x_2,s)}\cdot {\rm det}\dfrac{\partial(f_1,x_2,s)}{\partial(x_1,x_2,x_3)} = 1A(2/A) = -2$.  

 It is easy to see that 
$$ F(0, 0, -1/4) = F(1, -3/2, 13/2) = F(-1/3, 3/2, 13/2) = (-1/4, 0, 0).$$ 
 
\bigskip

 Here,  I do calculate its Jacobian determinant ${\rm det}(\partial f_j/\partial x_i)$  \uline{directly} using  calculation writing  (long-form calculation) with only Mathematica v.3  in the sequel.  
 I do not know why nobody does it. 

%\bigskip

-------------

$\bullet\ f_1[x_1,x_2,x_3]:=(1+x_1x_2)^3x_3+x_2^2(1+x_1x_2)(4+3x_1x_2).$

\bigskip

$f_{11}={\rm Expand}[{\rm D}[f_1[x_1,x_2,x_3],x_1]]
=7x_2^3 + 6x_1x_2^4 + 3x_2x_3 + 6x_1x_2^2x_3 + 3x_1^2x_2^3x_3.$

$f_{12}={\rm Expand}[{\rm D}[f_1[x_1,x_2,x_3],x_2]]
 = 8x_2 + 21x_1x_2^2 + 12x_1^2x_2^3 + 3x_1x_3 + 6x_1^2x_2x_3 + 3x_1^3x_2^2x_3.$

$f_{13}={\rm Expand}[{\rm D}[f_1[x_1,x_2,x_3],x_3]]
= 1 + 3x_1x_2 + 3x_1^2x_2^2 + x_1^3x_2^3.$

%\bigskip

-------------
\bigskip

$\bullet\  f_2[x_1,x_2,x_3]:=x_2+3x_1(1+x_1x_2)^2x_3+3x_1x_2^2(4+3x_1x_2).$

\bigskip

$f_{21}={\rm Expand}[{\rm D}[f_2[x_1,x_2,x_3],x_1]] = 12x_2^2 + 18x_1x_2^3 + 3x_3 + 12x_1x_2x_3 + 9x_1^2x_2^2x_3.$

$f_{22}={\rm Expand}[{\rm D}[f_2[x_1,x_2,x_3],x_2]] = 1 +24x_1x_2 + 27x_1^2x_2^2 + 6x_1^2x_3 + 6x_1^3x_2x_3.$

$f_{23}={\rm Expand}[{\rm D}[f_2[x_1,x_2,x_3],x_3]] = 3x_1 + 6x_1^2x_2 + 3x_1^3x_2^2.$ 

%\bigskip     
  
------------

\bigskip

$\bullet\ f_3[x_1,x_2,x_3]:=2x_1-3x_1^2x_2-x_1^3x_3.$

\bigskip

$f_{31}={\rm D}[2x_1-3x_1^2x_2-x_1^3 x_3,x_1] = 2 - 6x_1x_2 - 3x_1^2x_3.$

$f_{32}={\rm D}[2x_1-3x_1^2x_2-x_1^3x_3,x_2] = -3x_1^2.$

$f_{33}={\rm D}[2x_1-3x_1^2x_2-x_1^3x_3,x_3] = -x_1^3.$
  
%\bigskip

---------------------------------

\bigskip
$\bullet\ {\rm det Jac}(f_1,f_2,f_3;x_1,x_2,x_3)$\ \ \ \ (Jacobian determinant)

$= {\rm det}\left(\dfrac{\partial f_i}{\partial x_j}\right)$

$={\rm Det}[\{\{f_{11},f_{12},f_{13}\},\{f_{21},f_{22},f_{23}\},\{f_{31},f_{32},f_{33}\}\}]$

$= -f_{13}f_{22}f_{31} + f_{12}f_{23}f_{31}+f_{13}f_{21}f_{32}-f_{11}f_{23}f_{32}-f_{12}f_{21}f_{33}+f_{11}f_{22}f_{33}$

$- -2$.

%\bigskip

-----------------
 
\bigskip

More precisely, 

put

$j:= -f_{13}f_{22}f_{31} + f_{12}f_{23}f_{31} + f_{13}f_{21}f_{32} - f_{11}f_{23}f_{32} - f_{12}f_{21}f_{33} + f_{11}f_{22}f_{33}.$
 
%\bigskip

-------------------------------

\bigskip

$ j_1 = {\rm Expand}[-f_{13}f_{22}f_{31}]  = -2 -  48x_1x_2 - 42x_1^2 x_2^2 + 304x_1^3x_2^3 +  714x_1^4x_2^4 + 576x_1^5x_2^5 + 162x_1^6x_2^6 -  9x_1^2x_3 + 69x_1^3x_2x_3 + 378x_1^4x_2^2x_3 + 
      630 x_1^5 x_2^3 x_3 + 447 x_1^6 x_2^4 x_3 + 117 x_1^7 x_2^5 x_3 + 
      18 x_1^4 x_3^2 + 72 x_1^5 x_2 x_3^2 + 108 x_1^6 x_2^2 x_3^2 + 
      72 x_1^7 x_2^3 x_3^2 + 18 x_1^8 x_2^4 x_3^2.$

\bigskip
 
---

$j_2 = {\rm Expand}[f_{12}f_{23}f_{31}]
   = 48x_1 x_2 + 78x_1^2x_2^2 - 294x_1^3x_2^3 - 846x_1^4x_2^4 - 
    738x_1^5x_2^5 - 216x_1^6x_2^6 + 18x_1^2x_3 - 
    54x_1^3x_2x_3 - 441x_1^4x_2^2x_3 - 810x_1^5x_2^3x_3 - 
    603x_1^6x_2^4x_3 - 162x_1^7x_2^5x_3 - 27x_1^4x_3^2 - 
    108x_1^5x_2x_3^2 - 162x_1^6x_2^2x_3^2 - 
    108x_1^7x_2^3x_3^2 - 27x_1^8x_2^4x_3^2.$

\bigskip

----

$j_3={\rm Expand}[f_{13}f_{21}f_{32}]
   = -36 x_1^2 x_2^2 - 162 x_1^3 x_2^3 - 270 x_1^4 x_2^4 - 
      198 x_1^5 x_2^5 - 54 x_1^6 x_2^6 - 9 x_1^2 x_3 - 
      63 x_1^3 x_2 x_3 - 162 x_1^4 x_2^2 x_3 - 198 x_1^5 x_2^3 x_3 -117 x_1^6 x_2^4 x_3 - 27 x_1^7 x_2^5 x_3.$

\bigskip

----

$j_4 = {\rm Expand}[-f_{11}f_{23}f_{32}]
      = 63 x_1^3 x_2^3 + 180 x_1^4 x_2^4 + 171 x_1^5 x_2^5 + 
        54 x_1^6 x_2^6 + 27 x_1^3 x_2 x_3 + 108 x_1^4 x_2^2 x_3 + 
        162 x_1^5 x_2^3 x_3 + 108 x_1^6 x_2^4 x_3 + 
        27 x_1^7 x_2^5 x_3.$

\bigskip

----

$j_5 ={\rm Expand}[-f_{12} f_{21} f_{33}]
 = 252 x_1^4 x_2^4 + 522 x_1^5 x_2^5 + 216 x_1^6 x_2^6 + 
        99 x_1^4 x_2^2 x_3 + 414 x_1^5 x_2^3 x_3 + 
        477 x_1^6 x_2^4 x_3 + 162 x_1^7 x_2^5 x_3 + 9 x_1^4 x_3^2  + 54 x_1^5 x_2 x_3^2 + 108 x_1^6 x_2^2 x_3^2 + 
        90 x_1^7 x_2^3 x_3^2 + 27 x_1^8 x_2^4 x_3^2.$

\bigskip

----

$j_6 = {\rm Expand}[f_{11} f_{22} f_{33}]
 =-7 x_1^3 x_2^3 - 174 x_1^4 x_2^4 - 333 x_1^5 x_2^5 - 
        162 x_1^6 x_2^6 - 3 x_1^3 x_2 x_3 - 78 x_1^4 x_2^2 x_3 - 
        270 x_1^5 x_2^3 x_3 - 312 x_1^6 x_2^4 x_3 - 
        117 x_1^7 x_2^5 x_3 - 18 x_1^5 x_2 x_3^2 - 
        54 x_1^6 x_2^2 x_3^2 - 54 x_1^7 x_2^3 x_3^2 - 
        18 x_1^8 x_2^4 x_3^2.$

\bigskip

------------------------------------------- 

Then 

\bigskip

$j := j_1 + j_2 + j_3 + j_4 + j_5 + j_6$

$ =  -f_{13} f_{22} f_{31} + f_{12} f_{23} f_{31} + f_{13} f_{21} f_{32} - f_{11} f_{23} f_{32}
 -f_{12} f_{21} f_{33} + f_{11} f_{22} f_{33}$

 $= {\rm Simplify}[\\
(-2 -  48 x_1 x_2 - 42 x_1^2 x_2^2 + 304 x_1^3 x_2^3 + 714 x_1^4 x_2^4 + 576 x_1^5 x_2^5 + 162 x_1^6 x_2^6 - 9 x_1^2 x_3 + 69 x_1^3 x_2 x_3 + 378 x_1^4 x_2^2 x_3 + 630 x_1^5 x_2^3 x_3 + 447 x_1^6 x_2^4 x_3 + 117 x_1^7 x_2^5 x_3 + 18 x_1^4 x_3^2 + 72 x_1^5 x_2 x_3^2 + 108 x_1^6 x_2^2 x_3^2 + 72 x_1^7 x_2^3 x_3^2 + 18 x_1^8 x_2^4 x_3^2)\\
 + (48x_1 x_2 + 78x_1^2x_2^2 - 294x_1^3x_2^3 - 846x_1^4x_2^4 - 
    738x_1^5x_2^5 - 216x_1^6x_2^6 + 18x_1^2x_3 - 
    54x_1^3x_2x_3 - 441x_1^4x_2^2x_3 - 810x_1^5x_2^3x_3 - 
    603x_1^6x_2^4x_3 - 162x_1^7x_2^5x_3 - 27x_1^4x_3^2 - 
    108x_1^5x_2x_3^2 - 162x_1^6x_2^2x_3^2 - 
    108x_1^7x_2^3x_3^2 - 27x_1^8x_2^4x_3^2)\\
 + (-36 x_1^2 x_2^2 - 162 x_1^3 x_2^3 - 270 x_1^4 x_2^4 - 198 x_1^5 x_2^5 - 54 x_1^6 x_2^6 - 9 x_1^2 x_3 - 63 x_1^3 x_2 x_3 - 162 x_1^4 x_2^2 x_3 - 198 x_1^5 x_2^3 x_3 - 117 x_1^6 x_2^4\ x_3 - 27 x_1^7 x_2^5 x_3) \\
  + (63 x_1^3 x_2^3 + 180 x_1^4 x_2^4 + 171 x_1^5 x_2^5 + 54 x_1^6 x_2^6 +  27 x_1^3 x_2 x_3 + 108 x_1^4 x_2^2 x_3 + 162 x_1^5 x_2^3 x_3 + 108 x_1^6 x_2^4 x_3 + 27 x_1^7 x_2^5 x_3)\\
  + (252 x_1^4 x_2^4 + 522 x_1^5 x_2^5 + 216 x_1^6 x_2^6 + 99 x_1^4 x_2^2 x_3 + 414 x_1^5 x_2^3 x_3 + 477 x_1^6 x_2^4 x_3 + 162 x_1^7 x_2^5 x_3 + 9 x_1^4 x_3^2 + 54 x_1^5 x_2 x_3^2 + 108 x_1^6 x_2^2 x_3^2 + 90 x_1^7 x_2^3 x_3^2 + 27 x_1^8 x_2^4 x_3^2) \\
  + (-7 x_1^3 x_2^3 - 174 x_1^4 x_2^4 - 333 x_1^5\ x_2^5 - 162 x_1^6 x_2^6 - 3 x_1^3 x_2 x_3 - 78 x_1^4 x_2^2 x_3 - 270 x_1^5 x_2^3 x_3 - 312 x_1^6 x_2^4 x_3 - 117 x_1^7 x_2^5 x_3 - 18 x_1^5 x_2 x_3^2 - 54 x_1^6 x_2^2 x_3^2 - 54 x_1^7 x_2^3 x_3^2 - 18 x_1^8 x_2^4 x_3^2)] $ \\

%\bigskip

$=-2 \ \ (\in \mathbb{C}^\times)$.

%\bigskip

-----------------

 Therefore \uline{the polynomial map $(f_1,f_2,f_3) : \mathbb{C}^3 \rightarrow \mathbb{C}^3$ is  \'{e}tale}, that is,  a Keller map.

 %%%%%%%%%%%%%%%%%%%%

%\vspace{5mm}
%
%\noindent
%{\bf Acknowledgment :}  The author would like to be grateful to Moeko ODA for walking with him for a long time, to his son Shuhei ODA and to  young grandsons Naoki ODA and Takuma KUSHIDA, and finally to Unyo ODA, the beloved dog of his family, for having always cheered him up (who passed away on July 1, 2014 in Kochi City, JAPAN). 

 %%%%%%%%%%%%%%               

\vspace{7mm}

\begin{flushright}

\begin{minipage}{100mm}
\begin{small}

\noindent
``\begin{it}
 ${}^{6}$In the morning sow your seed, and at evening withhold not your 

\ \ hand,

\ \ for you do not know which will prosper, this or that, or whether both

\ \ alike will be good.
\end{it}''

\rightline{--------- ECCLESIASTES 11 \ (ESV)}   

\bigskip 
\bigskip

\noindent
`` \begin{it}
 ${}^{11}$Again I saw that under the sun 

 \ \ the race is nor to the swift,

 \ \ nor the battle to the strong,

 \ \ nor bread to the wise,

 \ \ nor riches to the intelligent,

 \ \ nor favour to those with knowledge,

 \ \ but time and chance happen to them all.

 ${}^{12}$For man does not know his time. $\cdots\cdots$  
 \end{it}
''

\rightline{--------- ECCLESIASTES 9 \ (ESV)}   

%\bigskip 
%\bigskip
%
%\noindent
%``\begin{it}
%${}^{14}$For He Himself knows our frame; 
%
%\begin{eqnarray}
%
%\end{eqnarray}
%
%\ \ \ \ He is mindful that we are (made of) but dust.
%
% \ \ ${}^{15}$As for man, his days are like grass;
%
%\ \ \ \ As a flower of the field, so he flourishes. 
%
%\ \  ${}^{16}$When the wind has passed over it,
%
%\ \ \ \ it is no more, and its place acknowledges it no longer.
%\end{it}
%''
%
%\rightline{ --------- PSALM 103, 14-16.\ (NASB)}
\end{small}
\end{minipage}
\end{flushright}
 
%%%%%%%%%%%%%%%%%%%%%%%%

\end{document}